\documentclass[11pt]{article}

\usepackage{amsmath}
\usepackage{amssymb}
\usepackage{amsfonts}
\usepackage{amsthm}
\usepackage{hyperref}
\usepackage{url}

\newcommand{\eqn}[2]{\begin{equation}\label{#1}#2\end{equation}}
\newcommand{\eqnst}[1]{\begin{equation*}#1\end{equation*}}
\newcommand{\eqnspl}[2]{\begin{equation}\begin{split}\label{#1}%
   #2\end{split}\end{equation}}
\newcommand{\eqnsplst}[1]{\begin{equation*}\begin{split}%
   #1\end{split}\end{equation*}}

\def\bone{\mathbf{1}}

\def\bp{\mathbf{p}}
\def\P{\mathbf{P}}
\def\Prob{\mathbf{P}}

\def\E{\mathbf{E}}
\def\cF{\mathcal{F}}
\def\cG{\mathcal{G}}

\theoremstyle{plain}
\newtheorem{theorem}{Theorem}[section]
\newtheorem{proposition}[theorem]{Proposition}
\newtheorem{lemma}[theorem]{Lemma}

\theoremstyle{definition}

\numberwithin{equation}{section}

\title{Asymptotics of the optimum in
discrete sequential assignment}
\author{Antal A.~J\'arai\thanks{Department of Mathematical Sciences, University of Bath,
Bath BA2 7AY, United Kingdom. E-mail: {\tt A.Jarai@bath.ac.uk}}\\ University of Bath}

\begin{document}

\maketitle

\begin{abstract}
We consider the stochastic sequential assignment problem of Derman, Lieberman 
and Ross (1972) corresponding to a discrete distribution supported on a finite
set of points. We use large deviation estimates to compute the asymptotics of 
the optimal policy as the number of tasks $N \to \infty$. 
\end{abstract}



\section{Introduction}
\label{sec:intro}

Let us start with a simple example of the question we are interested in. Consider the 
following game: you start with a row of $N$ empty boxes, and a fair dice is rolled 
$N$ times. After each roll in turn, you have to assign the rolled value to one of 
the boxes. After all the rolls have taken place, an $N$-digit number is obtained, 
that you are trying to maximize (see \cite{J16} for a popularizing account of this game).
Let us be more precise about the sense of maximizing, for which the following are
two natural possibilities: (i) maximize the score on average, that is, the expected 
value of the final $N$-digit number; (ii) maximize the probability of obtaining the 
`perfect score', that is, when the assignment is monotone increasing from right
to left. 

\begin{figure}[htpb]
\centerline{%
\setlength{\unitlength}{0.15cm}
\begin{picture}(100,6)
\put(0,1){\framebox(4,4){6}}
\put(4,1){\framebox(4,4){}}
\put(8,1){\framebox(4,4){}}
\put(12,1){\framebox(4,4){5}}
\put(16,1){\framebox(4,4){}}
\put(20,1){\framebox(4,4){}}
\put(24,1){\framebox(4,4){}}
\put(28,1){\framebox(4,4){4}}
\put(32,1){\framebox(4,4){4}}
\put(36,1){\framebox(4,4){}}
\put(40,1){\framebox(4,4){}}
\put(44,1){\framebox(4,4){}}
\put(48,1){\framebox(4,4){3}}
\put(52,1){\framebox(4,4){3}}
\put(56,1){\framebox(4,4){3}}
\put(60,1){\framebox(4,4){}}
\put(64,1){\framebox(4,4){}}
\put(68,1){\framebox(4,4){}}
\put(72,1){\framebox(4,4){}}
\put(76,1){\framebox(4,4){2}}
\put(80,1){\framebox(4,4){}}
\put(84,1){\framebox(4,4){}}
\put(88,1){\framebox(4,4){}}
\put(92,1){\framebox(4,4){1}}
\put(96,1){\framebox(4,4){1}}
\end{picture}
}
\caption{Illustration of the setup in which a phase transition occurs. Here $N = 25$, $N_0 = 10$, and there are
$15$ rounds to go. In the picture $(n_1, \dots, n_5) = (3, 4, 3, 3, 2)$.}
\label{fig:contiguous}
\end{figure}
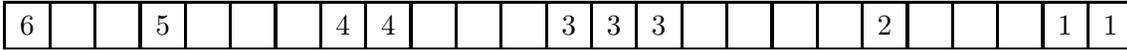

In the formulation (i), the problem is a special case of stochastic sequential 
assignment, introduced by Derman, Lieberman and Ross in \cite{DLR72}. They gave a 
recursive method to compute the optimal policy (when the values `rolled' are 
i.i.d.~samples from any distribution on $\mathbb{R}$ with finite mean). Formulation
(ii) is of interest in the discrete example outlined in the previous paragraph. 
An especially interesting feature of it is that it undergoes a controllability 
phase transition. In order to state this, a slightly different initial setup is convenient:
assume there are $N - N_0$ rounds to go, and the remaining set of empty boxes splits into 
$5$ contiguous intervals of lengths $n_1, \dots, n_5$, where $n_1 + \dots + n_5 = N-N_0$,
separated by contiguous intervals already filled by the values $2, \dots, 5$; 
see Figure \ref{fig:contiguous}. It follows from the main result of 
\cite{J17} that as $N-N_0 \to \infty$, if the vector $(n_i/(N-N_0))_{i=1}^5$ lies in the 
interior of a certain parallelepiped $K$, then the probability of obtaining the 
perfect score under the optimal policy is exponentially close to a constant $c > 0$,
whereas, if the same vector is in the exterior of $K$, it goes to $0$ 
exponentially. This critical phenomenon extends to a more general setting, where 
the $n_i$'s are assigned to edges of a finite graph, and the values `rolled' are the
vertices of the graph; see \cite{J17}. Estimates proved in \cite{J17} suggest that in
the critical region, when $(n_i/(N-N_0))_i$ is near $\partial K$, one can expect the
value of the optimal policy to exhibit Gaussian behaviour. Information on the asymptotics 
of the optimal policy in this region would help establish Gaussian behaviour rigorously.

A step towards the above goal is to determine the asymptotic behaviour, as 
$N \to \infty$, of the optimal policy in both problems (i) and (ii). This turns out to 
be easier for (i), and is the aim of the present paper. The reason why (i) is easier to 
analyze is the following. Let $Z_1, \dots, Z_N$ denote the values assigned to boxes 
$1, \dots, N$, respectively. Then in the case of (i), the objective is to maximize a 
\emph{linear} function of the $Z_i$'s (in expectation), wheras in the case of (ii) 
the objective is to maximize a non-linear function (the indicator of 
$Z_1 \le \dots \le Z_N$). We expect that in many, but not all discrete valued 
assignement problems, the formulations (i) and (ii) behave similarly, and we consider 
the case (i) as the first step in understanding (ii). 

Let us from now on restrict to the objective of maximizing the expected score, and ask:
when $N$ is large, where are the optimal locations to place the first rolled number,
if it is, respectively, $1, \dots, 6$? We will see that these are approximately
given by 
\eqnsplst
{ 1, \qquad 
    N \frac{\log \frac{4}{5}}{\log \frac{2}{5}}, \qquad
    N \frac{\log \frac{3}{4}}{\log \frac{1}{2}}, \qquad
    N \frac{\log \frac{2}{3}}{\log \frac{1}{2}}, \qquad
    N \frac{\log \frac{1}{2}}{\log \frac{2}{5}}, \qquad
    N, }
respectively, with boxes counted from the right end.

More generally, we consider the problem when the dice has $k \ge 3$ sides, with
real values $x_1 < \dots < x_k$ written on them, and the probability of rolling $x_i$ 
is $p_i > 0$, $i = 1, \dots, k$. 
As mentioned earlier, the optimal policy was found in \cite{DLR72}.
When the values to be assigned come from a \emph{continuous} distribution, it was shown in
\cite{AD72} by a law of large numbers argument, that the asymptotic profile of the 
optimal policy is given in terms of the scaled quantiles of the underlying distribution. 
We were surprised not to find in the literature any result on the asymptotics in the 
discrete case stated above. The discrete case is interesting, since 
the location of the optimum is determined by large deviation events, in contrast with 
the continuous case. 

The structure of the paper is the following.
In Section \ref{ssec:model} we precisely define the model considered, and state our
main result. Section \ref{ssec:optimum} gives the short large deviation computation 
that yields the asymptotic optimum values. The proof of the main theorem is given in 
Section \ref{sec:proof}.

\subsection{The model}
\label{ssec:model}

Let $\bp = (p_i)_{i=1}^{k}$ be a discrete distribution supported on the points
$x_1 < \dots < x_k$, such that $p_i > 0$ for all $1 \le i \le k$. 
Let $N \ge 1$, and let $r_N > \dots > r_1$ be given numbers (rewards).
Let $X_1, \dots, X_N$ be i.i.d.~random variables with distribution $\bp$.
Let $\cF_t = \sigma\{ X_1, \dots, X_t \}$, $0 \le t \le N$.
By a \textbf{policy} we mean a sequence of random indices 
$J(1), \dots, J(N) \in \{ 1, \dots, N \}$ such that:\\
(i) $J(t)$ is $\cF_t$-measurable; and\\
(ii) $J(1), \dots, J(N)$ are distinct.\\
By the \textbf{reward} of the policy we mean the random quantity
\eqnst
{ R(X;J)
  = \sum_{t=1}^N X_t r_{J(t)}. }
The aim is to maximize the expected reward. 
The theorem below follows directly from a theorem proved by 
Derman, Lieberman and Ross \cite{DLR72}; see also 
\cite[Section VI.7]{Ross-book}.
We denote by $1 = \ell_N(1) \le \ell_N(2) \le \dots \le \ell_N(k-1) \le \ell_N(k) = N$ 
the optimal locations for placing the values $x_i$, $i = 1, \dots, k$, 
among $N$ remaining empty spaces, $N \ge 1$ (making an arbitrary choice in case 
the optimal policy has ties).

\begin{theorem}{\cite{DLR72}}
\label{thm:DLR}
There exist numbers
\eqn{e:bound-ineq}
{ -\infty 
  = a_{N,0} 
  < \dots 
  < a_{N,N}
  = \infty, \quad N \ge 1, }
that only depend on $\bp$ and $(x_i)_{i=1}^k$ (and not on the $r_i$'s),
such that any optimal policy satisfies
\eqn{e:opt-ineq}
{ a_{N,\ell_{N}(i)-1} 
  \le x_i 
  \le a_{N,\ell_N(i)}, \quad 1 \le i \le k,\, N \ge 1. }
Moreover, we can choose 
\eqn{e:opt-choice}
{ a_{N,n}
  = \E [ X_{T(n)} ], \quad 1 \le n \le N-1,\, N \ge 2, } 
where $T(n)$ is the 
unique index $1 \le t \le N-1$ such that $J(t) = n$, and
$J$ is any optimal policy with $N-1$ rounds to go.
\end{theorem}

We are ready to formulate our main result. 
Let $f_i := \sum_{j \le i} p_j$, $i = 0, \dots, k$, and let
\eqn{e:ai}
{ d_i
  := \frac{ \log \frac{1 - f_i}{1 - f_{i-1}} }
     { \log \frac{f_{i-1} \, (1 - f_i)}{(1 - f_{i-1}) \, f_i} }, 
     \quad i = 2, \dots, k-1. }
We extend this definition to $i = 1$ and $i = k$ by setting 
$d_1 = 0$ and $d_k = 1$.

\begin{theorem}
\label{thm:main}
Let $k \ge 3$, and let $\bp$ be a positive probability vector.
Then any optimal policy satisfies
\eqnst
{ \ell_N(i)
  = N (d_i + o(1)), \quad \text{as $N \to \infty$, $i = 1, \dots, k$.} }
\end{theorem}

\subsection{Computation of the optimum}
\label{ssec:optimum}

The asymptotic optimum $d_i$ arises from the equality of two 
large deviation rates. We give this short computation before 
delving into the proof of Theorem \ref{thm:main}.
For $i = 2, \dots, k-1$ let us look for the value of 
$y \in (f_{i-1},f_{i})$ that yields
equality of the large deviation rates:
\eqnsplst
{ I^-_i(y) 
  &:= - \lim_{N \to \infty} \frac{1}{N} 
     \log \Prob \left[ \sum_{t=1}^N \bone_{X_t < x_i} \ge y N \right] \\
  I^+_i(y) 
  &:= - \lim_{N \to \infty} \frac{1}{N} 
     \log \Prob \left[ \sum_{t=1}^N \bone_{X_t > x_i} \ge (1-y) N \right]. }
Using the large deviation rate function for binomial random variables \cite{dH-book}, 
we have
\eqnsplst
{ I^-_i(y)
  &= y \log y + (1-y) \log (1-y) 
    - y \log f_{i-1} - (1-y) \log (1-f_{i-1}) \\
  I^+_i(y)
  &= y \log y + (1-y) \log (1-y) 
    - y \log f_{i} - (1-y) \log (1-f_{i}). } 
Therefore, equality occurs for $y = d_i$.

\section{Proof of the main theorem}
\label{sec:proof}

In this section we give the proof of Theorem \ref{thm:main}.
The idea of the proof is as follows. In order to arrive at a
contradiction, assume that $\limsup_{N \to \infty} \ell_N(i)/N > d_i$,
and $i$ is smallest with this property.
Let $U_N$ denote the random value that is assigned to the location
$\ell_N(i)-1$, that is, $U_N = X_{T(\ell_N(i)-1)}$. 
If we can deduce that $\E [ U_N ] > x_i$ for some 
(sufficiently large) $N$, this contradicts
Theorem \ref{thm:DLR}. Let $\tau$ be the random time when location 
$\ell_N(i)-1$ is filled, which is a stopping time:
\eqnst
{ \tau
  = \inf \{ t \ge 1 : J(t) = \ell_N(i)-1 \}. }
Let $I_t$ be the index of location $\ell_N(i)-1$ among the 
`empty spaces' at time $t$, that is:
\eqnst
{ I_t
  = \# \{ 1 \le n \le \ell_N(i)-1 : \text{$n \not= J(s)$ for any $1 \le s \le t$} \},
    \qquad 0 \le t < \tau. } 
Let $\sigma$ be the stopping time
\eqnst
{ \sigma
  = \inf \{ 0 \le t < \tau : I_t = \ell_{N-t}(i)+2 \} \in \{ 0, \dots, \tau-1 \} \cup \{ \infty \}. }
Observe that by definition, either $\sigma < \tau$ or $\sigma > \tau$ (when the set of $t$'s 
considered in the $\inf$ is empty).

We prove the following estimates.

\begin{proposition}
\label{prop:E-}
Let $2 \le i \le k-1$, and assume that:\\
(i) $\limsup_{N \to \infty} \ell_N(j)/N \le d_j$ for $1 \le j \le i-1$;\\
(ii) $\limsup_{N \to \infty} \ell_N(i)/N > b > d_i$ for some $b$.\\
Then there exist $b'$ satisfying $b > b' > d_i$ such that 
for all $N$ such that $\ell_N(i) \ge N b + 1$, we have 
\eqnst
{ \P [ U_N \le x_{i-1} ]
  \le \exp \left( - (1 + o(1)) \, I_i^- \left( b' \right) N \right), \qquad
      \text{as $N \to \infty$.} }
\end{proposition}

\begin{proposition}
\label{prop:E+}
Let $2 \le i \le k-1$, and assumme that:\\
(i) $\limsup_{N \to \infty} \ell_N(i)/N > b > d_i$ for some $b$.\\
Then there exists $b''$ satisfying $b > b'' > d_i$ and a constant 
$c_1 = c_1 ( b, x_i, x_{i+1} ) > 0$
such that for an infinite set $\cG$ of positive integers we have\\
(a) $\ell_N(i) \ge N b + 1$ for all $N \in \cG$; \\
(b) 
\eqnspl{e:Eplus-bound}
{ \E [ U_N \,|\, \sigma < \tau ]
  \ge x_i + c_1 \exp ( - (1 + o(1)) \, I^+_i \left( b'' \right) N ), \quad
      \text{for all $N \in \cG$.} }
(c) Moreover, there exists a constant $c_2 > 0$ such that 
$\P [ \sigma < \tau ] \ge c_2$ for all $N \in \cG$.
\end{proposition}

Let us first show that the two propositions imply the main theorem.

\begin{proof}[Proof of Theorem \ref{thm:main}.]
Adding a constant to the random variables $X_i$ changes the reward of any policy 
by a constant, and hence we may assume without loss of generality that 
$0 = x_1 < x_2 < \dots < x_k$.

It is sufficient to show that $\limsup_{N \to \infty} \ell_N(i)/N \le d_i$,
for $1 \le i \le k$, since then the analogous inequality for the liminf
follows by symmetry. In order to arrive at a contradiction, assume that $i$
is the smallest index such that the inequality fails. Then $2 \le i \le k-1$,
and we can fix a number $b$ with $\limsup_{N \to \infty} \ell_N(i)/N > b > d_i$,
so that Propositions \ref{prop:E-} and \ref{prop:E+} can be applied. Consider 
the subsequence $\cG$ provided by Proposition 2.2, and observe that due to 
part (a) of this proposition, the conclusion of Proposition 2.1 also holds for all
$N \in \cG$.

Since $\sigma \not= \tau$, we have 
\eqnsplst
{ \E [ U_N ]
  = \E \left[ U_N \, \bone_{\tau < \sigma} \right]
    + \E \left[ U_N \, \bone_{\sigma < \tau} \right], \quad N \in \cG. }
Proposition \ref{prop:E-} implies that 
\eqnsplst
{ \P [ U_N \le x_{i-1},\, \tau < \sigma ]
  \le \P [ U_N \le x_{i-1} ]
  \le \exp \left( -(1 + o(1)) \, I^-_i(b') \, N \right), \quad N \in \cG. }
Writing $\alpha = \P [ \sigma < \tau ]$ this implies
\eqnspl{e:lb-1}
{ \E \left[ U_N \, \bone_{\tau < \sigma} \right]
  &\ge x_i \P [ U_N \ge x_i,\, \tau < \sigma ] \\
  &\ge x_i \, (1 - \alpha) - x_i \, \exp \left( - (1 + o(1)) \, I^-_i(b') \, N \right), 
	  \quad N \in \cG. }
On the other hand, Proposition \ref{prop:E+} implies that 
\eqnspl{e:lb-2}
{ \E \left[ U_N \, \bone_{\sigma < \tau} \right]
  = \alpha \E [ U_N \,|\, \sigma < \tau ]
  \ge \alpha \left( x_i + c_1 \exp \left( -(1 + o(1)) \, I^+_i(b'') \, N \right) \right), \quad N \in \cG. } 
Putting together \eqref{e:lb-1} and \eqref{e:lb-2}, and using that 
$I^+_i(b'') < I^+_i(d_i) = I^-_i(d_i) < I^-_i(b')$, we get that for a
sufficiently large $N \in \cG$ we have
\eqnst
{ \E [ U_N ] > x_i. }
This is the desired contradiction, and completes the proof.
\end{proof}

Our next task is to prove Proposition \ref{prop:E-}. 
For what follows, recall the definition of $I_t$ from the beginning 
of this section, and observe that we have $I_0 = \ell_N(i)-1$, 
and the time evolution of $(I_t)_{0 \le t < \tau}$ is given by:
\eqn{e:I_t-evolution}
{ I_{t+1}
  = \begin{cases}
    I_t     & \text{when $X_{t+1} = x_j$ and $I_t < \ell_{N-t}(j)$, $2 \le j \le k$;} \\
    I_t - 1 & \text{when $X_{t+1} = x_j$ and $I_t > \ell_{N-t}(j)$, $1 \le j \le k-1$;}
    \end{cases} }
moreover, we have
\eqnst
{ \tau
  = t+1, \quad \text{when $X_{t+1} = x_j$ and 
     $I_t = \ell_{N-t}(j)$, $1 \le j \le k$.} }

\begin{proof}[Proof of Proposition \ref{prop:E-}.]
We fix $\delta > 0$ such that 
\eqnst
{ \delta < \min \left\{ b - d_i,\, d_i - d_{i-1},\, f_{i-1} - d_{i-1} \right\}. }
Due to assumption (i), there exists $N_2 \ge 1$ such that 
for all $N \ge N_2$ we have $\ell_N(i-1)/N \le d_{i-1} + \delta$.
Consider the stopping time
\eqnsplst
{ \sigma_1
  &:= \inf \{ 0 \le t < \tau : I_t \le \ell_{N-t}(i-1) \}. }
When $N-t \ge N_2$, we have
\eqnspl{e:o-estimate1}
{ (\ell_N(i) - 1) - \ell_{N-t}(i-1)
  \ge N b - (N-t) \left( d_{i-1} + \delta \right). }
Let 
\eqnsplst
{ N' 
  &:= N \frac{ b - d_{i-1} - \delta }{ 1 - d_{i-1} - \delta }. }
A simple computation shows that when $t < N'$, 
the right hand side of \eqref{e:o-estimate1} is 
larger than $t$. Since $I_t$ can decrease by at most
$1$ at each time step, this implies the deterministic
inequality $\sigma_1 \ge N'$. 

When $N-t < N_2$ and $N \ge N_3 := 2 N_2 / (b - d_i)$, we have
\eqnspl{e:o-estimate2}
{ (\ell_N(i)-1) - \ell_{N-t}(i-1)
  \ge N b - N_2
  \ge N \frac{b+d_i}{2}. } 
Putting $b' := (b+d_i)/2$ and 
\eqnsplst
{ K_t
  := \begin{cases}
     N b - (N-t) \left( d_{i-1} + \delta \right) & \text{when $N' \le t \le N - N_2$;} \\
     N b'                                        & \text{when $N - N_2 < t \le N$,}
     \end{cases} }
the estimates \eqref{e:o-estimate1} and \eqref{e:o-estimate2} imply that 
for $N \ge N_3$ we have
\eqnspl{e:drift-bound}
{ \Prob [ U_N \le x_{i-1} ]
  &\le \Prob [ \sigma_1 < \tau ] 
  \le \Prob \left[ \text{$\exists$ $N' \le t < \tau$ such that 
      $I_t = I_0 - K_t$} \right]. }
Taking into account the time-evolution of $(I_t)_{0 \le t < \tau}$,
the right hand side of \eqref{e:drift-bound} is at most:
\eqnspl{e:sum-prob}
{ \Prob \left[ \text{$\exists$ $N' \le t \le N$ such that 
      $\sum_{s=1}^t \bone_{X_s \le x_{i-1}} = K_t$} \right]
  \le \sum_{t = N'}^N \Prob \left[ \sum_{s = 1}^t \bone_{X_s \le x_{i-1}} = K_t \right]. }
Let us write $t = x N$, and put $g(x) = K_t/t$. Observe that we have
$(b - d_{i-1} - \delta)/(1 - d_{i-1} - \delta) \le x \le 1$, and we have
\eqnst
{ g(x)
  = \begin{cases}
    \frac{K_t}{t} = \frac{b}{x} - \frac{1-x}{x} \left( d_{i-1} + \delta \right) 
       & \text{when $N' \le t \le N - N_2$;} \\
    \frac{K_t}{t} \ge \frac{b'}{x} 
       & \text{when $N - N_2 < t \le N$.}
    \end{cases} }
Using the large deviation rate function for binomial random variables \cite{dH-book}, 
we have
\eqnspl{e:rate-function}
{ \log \Prob \left[ \sum_{s = 1}^t \bone_{X_s \le x_{i-1}} = K_t \right]
  &\le N \left[ o(1) - x g(x) \log \frac{g(x)}{f_{i-1}}
      - x (1 - g(x)) \log \frac{1 - g(x)}{1 - f_{i-1}} \right] \\
  &=: N \left[ o(1) + x F(x) \right], 
     \quad \text{as $N \to \infty$.} }
Note that the $o(1)$ term is uniform in $x$, since $t \ge N'$, and 
$N'$ grows (linearly) with $N$.

We first show that the expression $x F(x)$ inside square brackets in the right hand side
of \eqref{e:rate-function} is increasing when $N' \le t \le N - N_2$. We have 
\eqnsplst
{ g'(x)
  &= - \frac{b - d_{i-1} - \delta}{x^2} 
  = - \frac{g(x)}{x} + \frac{d_{i-1} + \delta}{x} \\
  -x g'(x)
  &= g(x) - d_{i-1} - \delta. }
We also have
\eqnspl{e:F'-formula}
{ &\frac{d}{dx} \left[ x F(x) \right]
  = F(x) + x F'(x) \\
  &\quad = F(x) + x \left[ - g'(x) \log \frac{g(x)}{f_{i-1}}
     + g'(x) \log \frac{1 - g(x)}{1 - f_{i-1}} 
     - g(x) \frac{g'(x)}{g(x)} 
     - (1 - g(x)) \frac{-g'(x)}{1 - g(x)} \right] \\
  &\quad = F(x) + \left( g(x) - d_{i-1} - \delta \right) \log \frac{g(x)}{f_{i-1}} 
    + \left( (1 - g(x)) - (1 - d_{i-1} - \delta) \right) \log \frac{1-g(x)}{1-f_{i-1}} \\
  &\quad = - (d_{i-1} + \delta) \log \frac{g(x)}{f_{i-1}}
    - (1 - d_{i-1} - \delta) \log \frac{1-g(x)}{1-f_{i-1}}. }
Note that $b \le g(x) \le 1$. We show that with $c = d_{i-1} + \delta$ and
$d = f_{i-1}$, the function
\eqnst
{ h(y)
  = - c \log \frac{y}{d}
    - (1 - c) \log \frac{1-y}{1-d} }
is everywhere positive on the interval $y \in [b,1]$.
Indeed, we have $y \ge b > d_i > f_{i-1} = d > d_{i-1} + \delta = c$,
and 
\eqnst
{ h'(y)
  = - \frac{c}{y} + \frac{1 - c}{1-y}
  = \frac{y ( 1 - c ) - c (1 - y)}{y \, (1 - y)}
  = \frac{y - c}{y \, (1 - y)}
  > 0, \quad \text{ for $y \ge d$.} }
Therefore, $h(y) > h(d) = 0$.  
This gives that 
\eqn{e:main-bit}
{ \sup \{ x F(x) : N'/N \le x \le 1 - N_2/N \}
  \le -b \log \frac{b}{f_{i-1}} - (1-b) \log \frac{1-b}{1-f_{i-1}}
  = - I^-_i(b). }
On the other hand, we have
\eqn{e:side-bit}
{ \sup \{ x F(x) : 1 - N_2/N \le x \le 1 \}
  = - (1 + o(1)) I^-_i \left( b' \right), 
    \qquad \text{as $N \to \infty$.} }
The estimates \eqref{e:main-bit} and \eqref{e:side-bit} ensure that the 
right hand side of \eqref{e:sum-prob} is at most
\eqnst
{ \exp ( - (1 + o(1)) \, I^-_i( b' ) \, N ), \qquad \text{as $N \to \infty$.} }
This completes the proof.
\end{proof}

In order to complement the bound provided by Proposition \ref{prop:E-},
we seek a lower bound on the probability that a value $\ge i+1$ is 
assigned to $\ell_N(i)-1$.
To start, the following lemma provides a `continuity' result for the optimal policy.

\begin{lemma}
\label{lem:opt-continuity}
For every $N \ge 1$ we have $\ell_N(i) \in \{ \ell_{N-1}(i),\, \ell_{N-1}(i) - 1 \}$.
\end{lemma}

\begin{proof}
It is easy to verify from the proof of Theorem \ref{thm:DLR}
using induction on $N$, 
that for all $N \ge 2$ the numbers $a_{N,r}$, $r = 0, \dots, N$ are
strictly increasing. 

Fix $2 \le i \le k-1$.
Let $r \le \ell_{N-1}(i)-2$. Then due to the proof of 
Theorem \ref{thm:DLR} \cite{DLR72},
we have that $a_{N,r}$ is a convex combination of the
numbers $a_{{N-1},r-1} < a_{{N-1},r}$, both of which are 
at most $x_i$, and there is positive weight on the smaller value.
It follows that $a_{N,r} < x_i$, and hence $\ell_N(i) \ge \ell_{N-1}(i)-1$.

By a similar argument, if $r \ge \ell_{N-1}+1$, then 
$a_{N,r}$ is a convex combination of the numbers 
$a_{{N-1},r-1} < a_{{N-1},r}$, both of which are
at least $x_i$, and there is positive weight on the larger
value. Therefore, $a_{N,r} > x_i$ and $\ell_N(i) \le \ell_{N-1}(i)$.
\end{proof}

We will need to restrict $N$ to a set of `good' values.
We need the following lemma to define these.

\begin{lemma}
\label{lem:o-bound}
Let $2 \le i \le k-1$, let $d_i < v < b < 1$, and 
assume that $\limsup_{N \to \infty} \ell_N(i)/N > b$. Then 
there are infinitely many values of $N$ such that the 
following hold:\\
(i) $\ell_N(i) \ge b N + 1$; \\
(ii) we have
\eqn{e:o-bound}
{ \ell_{N-t}(i)
  \le \ell_{N}(i) - v t, \quad
      0 \le t \le N-1. }
\end{lemma}

\begin{proof}
There are infinitely many numbers $N'$ such that $\ell_{N'}(i) \ge b N' + 1$.
We first claim that with $c_1 = (b-v)/(1-v)$, and all $N'$ sufficiently large,
there exists $N$ with $c_1 N' \le N \le N'$ such that 
\eqn{e:o-bound1}
{ \ell_{N-t}(i)
  \le \ell_{N}(i) - v t, \quad
      0 \le t \le N-1. }
Should such $N$ not exist, there would be an integer $r \ge 1$ and a sequence of 
numbers $N' = n_0 > n_1 > \dots > n_r \ge 1$ with $1 \le n_r < c_1 N'$,
such that 
\eqnst
{ \ell_{n_s}(i) 
  > \ell_{n_{s-1}}(i) - v (n_{s-1} - n_{s}), \quad
    s = 1, \dots, r. }
(Indeed, such a sequence can be choosen inductively, using the negation of
\eqref{e:o-bound1} for $N = n_0, n_1, \dots$ in turn.)
This implies that 
\eqnsplst
{ n_r 
  &\ge \ell_{n_r}(i)
  = \ell_{N'}(i) + \sum_{s=1}^r (\ell_{n_s} - \ell_{n_{s-1}})
  > \ell_{N'}(i) - v \, \sum_{s=1}^r (n_{s-1} - n_s)
  = \ell_{N'}(i) - v (N' - n_r) \\
  &> b N' - v (N' - n_r)
  = v n_r + (b - v) N'. }
Hence $n_r \ge N' \, (b-v)/(1-v) = c_1 N'$.
This is a contradiction, and hence $N$ exists with 
the claimed property \eqref{e:o-bound1}.

Let now $N$ be the largest integer $\le N'$ with the 
property \eqref{e:o-bound1}. Then it also holds that
\eqn{e:o-bound2}
{ \ell_{N+t}(i) 
  \le \ell_{N}(i) + v t, \quad
      0 \le t \le N' - N. } 
Indeed, \eqref{e:o-bound2} holds automatically for $t = 0$, and should it 
be violated for some $1 \le t \le N' - N$, then with the smallest such  
$t$ the value $N + t$ would also satisfy \eqref{e:o-bound1}, 
contradicting the maximality of $N$.
Observe that due to \eqref{e:o-bound2} we have
\eqnst
{ \ell_{N}(i)
  \ge \ell_{N'}(i) - v (N'-N)
  \ge b N' + 1 - v (N'-N)
  \ge b N' + 1 - b (N'-N)
  = b N + 1, } 
and hence $N$ satisfies both properties in the Lemma.
\end{proof}

We will write $\cG = \cG(b,v,i)$ for the set of values of $N \ge 1$
such that both (i)--(ii) in Lemma \ref{lem:o-bound} hold.

\begin{lemma}
There exists $c_2 > 0$ such that for $N \in \cG$ we have
$\P [ \sigma < \tau ] \ge c_2$. 
\end{lemma}

\begin{proof}
Let $N \in \cG$. Due to Lemma \ref{lem:opt-continuity} and \eqref{e:o-bound}, 
there exists a deterministic $t_0 = t_0(N)$ such that on the event
$\{ X_1 = \dots = X_{t_0} = x_k \}$, we have $I_{t_0} = \ell_{N-t_0}(i)$.
Let us choose the smallest such $t_0$. Then there exists a deterministic
upper bound $T_0$, so that $t_0(N) \le T_0$
Due to Lemma \ref{lem:opt-continuity} and \eqref{e:o-bound},
There exists $t_1 = t_1(N) > t_0$ such that on the event
$\{ X_1 = \dots = X_{t_1} = x_k \}$,
we have $\sigma = t_1 < \tau$. Since there is also a deterministic upper
bound $t_1(N) \le T_1$, we have $\P [ \sigma < \tau ] \ge p_k^{T_1}$.
\end{proof}

We introduce the martingale $Y_t := \E [ U_N \,|\, \cF_t ]$.
Observe that $Y_0 = \E [ U_N ]$ and $Y_t = U_N$ when $t \ge \tau$.

\begin{proof}[Proof of Proposition \ref{prop:E+}.] 
On the event $\sigma < \tau$, consider the stopping time
\eqnst
{ \sigma''
  = \inf \{ t > \sigma : I_t = \ell_{N-t}(i) + 1 \}. }
On the event $\sigma < \sigma'' < \tau$, we have
\eqnst
{ \E [ U_N \,|\, \cF_{\sigma''} ]
  = Y(\sigma'') 
  \ge x_i. } 
On the event $\sigma < \tau \le \sigma''$, we have
$Y(\tau) \ge x_{i+1}$. We construct an event 
$E^+ \subset \{ \sigma < \tau \le \sigma'' \}
\subset \{ \sigma < \tau,\, U_N > x_i \}$
such that  
\eqn{e:Eplus-event}
{ \Prob [ E^+ ]
  \ge c_1 \exp ( - (1 + o(1))\, I^+_i(b'') \, N ). }
This will prove the proposition, since using the martingale
property we have
\eqnsplst
{ \E [ U_N \,|\, \sigma < \tau ]
  &= \frac{1}{\P [ \sigma < \tau ]}
    \left( \E \left[ U_N \, \bone_{\sigma < \sigma'' < \tau} \right]
    + \E \left[ U_N \, \bone_{\sigma < \tau \le \sigma''} \right] \right) \\
  &= \frac{1}{\P [ \sigma < \tau ]}
    \left( \E \left[ Y(\sigma'') \, \bone_{\sigma < \sigma'' < \tau} \right]
    + \E \left[ Y(\tau) \, \bone_{\sigma < \tau \le \sigma''} \right] \right) \\
  &\ge \frac{1}{\P [ \sigma < \tau ]}
    \left( x_i \, \P [ \sigma < \sigma'' < \tau ] 
    + x_{i+1} \, \P [ \sigma < \tau \le \sigma'' ] \right) \\ 
  &\ge x_i + \left( x_{i+1} - x_i \right) \, \frac{\Prob [ E^+ ]}{\P [ \sigma < \tau ]}
  \ge x_i + \left( x_{i+1} - x_i \right) \, \P [ E^+ ]. }

In order to define $E^+$ and prove \eqref{e:Eplus-event}, 
we fix a positive integer $A$ whose value will be determined
later. Due to \eqref{e:o-bound}, there exists an integer $t_2$
such that 
\eqnst
{ I_0 
  = \ell_{N} - v t_2 + A. }
Let 
\eqnsplst
{ E^+_1 
  &= \{ V_{1} = \dots = V_{t_1} = x_k \} \\
  E^+_2
  &= \{ V_{t_1+1} = \dots = V_{t_2} = x_k \} \\
  E^+
  &= E^+_1 \cap E^+_2 \cap \{ \tau \le \sigma'' \}. }
Note that on $E^+_1 \cap E^+_2$ we have 
$I_{t_2} = I_{t_1} = I_0 = \ell_{N} - v t_2 + A$.
Therefore, due to \eqref{e:o-bound} for $t > t_2$ we have
\eqnst
{ \ell_{N-t} 
  \le \ell_N - v t
  = \ell_N - v t_2 + A - (t - t_2) v - A
  = I_{t_2} - v (t - t_2) - A. }
This implies the inclusion of events:
\eqnsplst
{ E^+
  &\supset E^+_1 \cap E^+_2 \cap \left\{ \text{$\forall\, t \ge t_2$ we have 
    $\ell_{N-t}+1 < I_t$} \right\} \\
  &\supset E^+_1 \cap E^+_2 \cap \left\{ \text{$\forall\, t \ge t_2$ we have
    $I_{t_2} - v (t - t_2) - (A-1) < I_t$} \right\} \\
  &\supset E^+_1 \cap E^+_2 \cap \left\{ \text{$\forall\, t \ge t_2$ we have 
    $\sum_{t_2 < s \le t} \mathbf{1}_{X_s \le x_i} < A-1 + v (t - t_2)$} \right\}. }
Writing $Y_s = \mathbf{1}_{X_s \le x_i}$ and 
$G(t) = \{ \sum_{t_2 < s \le t} Y_s < A-1 + v (t - t_2) \}$
this gives
\eqnst
{ \P [ E^+ ]
  \ge p_k^{t_2} \, \P \left[ \sum_{t_2 < s \le t} Y_s < A-1 + v (t - t_2) 
      \, \text{for all $t \ge t_2$} \right]
  = p_k^{t_2} \P \left[ \cap_{t_2 < t \le N} G(t) \right]. }
Let us fix a number $b''$ that satisfies $d_i < b'' < v$, and 
put $M = \lfloor b'' (N - t_2) \rfloor$. 
Let $H = \{ \sum_{t_2 < s \le N} Y_s = M \} \subset G(N)$.
Then due to the definition of $I^+_i$, 
for sufficiently large $N$ we have
\eqnst
{ \P [ H ]  
  = \exp \left( - (1 + o(1)) \, I^+_i(b'') \, (N - t_2) \right)
  \ge \exp ( - (1 + o(1)) \, I^+_i(b'') \, N ). }
We show that 
\eqn{e:lb-required}
{ \P [ \cup_{t_2 < t < N} G(t)^c \,|\, H ]
  \le \sum_{t_2 < t < N} \P [ G(t)^c \,|\, H ]
  \le 1/2, }
if $A$ is large enough, which implies the required statement.
Under the conditioning on $H$, $Y_{t_2+1}, \dots, Y_N$ have the 
same law as an i.i.d.~sequence $Y'_{t_2+1}, \dots, Y'_N$ with
\eqnst
{ \P [ Y'_i = 1 ] 
  = \frac{M}{N - t_2} \qquad
  \P [ Y'_i = 0 ]
  = 1 - \frac{M}{N - t_2} }
conditioned on $\sum_{t_2 < s \le N} Y'_s = M$. 
Write $S'_t = \sum_{t_2 < s \le t} Y'_s$.
Then we need to give an upper bound on
\eqnst
{ \P [ S'_t \ge A-1 + v (t - t_2) \,|\, S'_N = M ], }

We distinguish the cases $N/2 \le t < N$ and 
$t_2 < t < N/2$. Since the mean of $Y'_s$ is less than $v$, 
when $N/2 \le t < N$ we have
\eqnst
{ \P [ S'_t \ge A-1 + v (t - t_2) \,|\, S'_N = M ]
  \le \frac{\P [ S'_t \ge A-1 + v (t - t_2) ]}{\P [ S'_N = M ]}
  \le C \sqrt{N} \, C \exp ( - c N ), }
where we used the local limit theorem to lower bound the
denominator. The right hand side sums to less than
$1/4$, if $N$ is large enough.
 
When $t_2 < t < N/2$, 
we can write
\eqnsplst
{ &\P [ S'_t \ge A-1 + v (t - t_2) \,|\, S'_N = M ]
  = \frac{\P [ S'_t \ge A-1 + v (t - t_2),\, S'_N = M ]}{\P [ S'_N = M ]} \\
	&\qquad = \frac{1}{\P [ S'_N = M ]} \sum_{y \ge A-1 + v (t - t_2)} \P [ S'_t = y,\, S'_N = M ] \\
	&\qquad = \frac{1}{\P [ S'_N = M ]} \sum_{y \ge A-1 + v (t - t_2)} \P [ S'_t = y ] \, 
	  \P [ S'_N = M \,|\, S'_t = y ]. }
For each fixed $y$, the last conditional probability is a binomial probability with $N-t$ trials, 
where the probability of a success is $\P [ Y'_i = 1 ] = \frac{M}{N - t_2}$. This success
probability is bounded away from $0$ and $1$, due to $M = \lfloor b'' (N - t_2) \rfloor$
and $0 < b'' < 1$. Hence $\P [ S'_N = M \,|\, S'_t = y ]$ is bounded above by the 
maximum of the above binomial distribution, which is $\le C / \sqrt{N - t} \le C / \sqrt{N/2}$,
uniformly in $y$. This gives that 
\eqnsplst
{ &\P [ S'_t \ge A-1 + v (t - t_2) \,|\, S'_N = M ]
  \le \frac{1}{\P [ S'_N = M ]} \sum_{y \ge A-1 + v (t - t_2)} \P [ S'_t = y ] \, 
	  \frac{C}{\sqrt{N}} \\
  &\qquad \le \P [ S'_t \ge A-1 + v (t - t_2) ] \, \frac{C / \sqrt{N}}{\P [ S'_N = M ]} 
  \le C \P [ S'_t \ge A-1 + v (t - t_2) ]. }
Since again the mean of $Y'_s$ is less than $v$, the right hand side is
summable in $t$. Moreover, by choosing $A$ large, we can make it
sum to a value less than $1/4$. The two cases together 
yield the required \eqref{e:lb-required}, and establish \eqref{e:Eplus-event}.
\end{proof}

\medbreak

\textbf{Acknowledgements.} Research supported by the University of Bath.

\end{document}